\documentclass[11pt]{article}

\usepackage{amsmath,amsthm,amsfonts}

\usepackage{amssymb}

\usepackage{epsfig}

\usepackage{rotating}

\usepackage{subfigure}

\usepackage{authblk}
\usepackage{color}

\begin{document}

\title{Improving solution accuracy and convergence for stochastic physics parameterizations with colored noise}

\author[1]{Panos Stinis} 
\author[3]{Huan Lei}
\author[1]{Jing Li} 
\author[2]{Hui Wan} 

\affil[1]{Advanced Computing, Mathematics and Data Division, Pacific Northwest National Laboratory, Richland WA 99354}
\affil[2]{Atmospheric Sciences and Global Change Division, Pacific Northwest National Laboratory, Richland WA 99354}
\affil[3]{Department of Computational Mathematics, Science and Engineering and Department of Statistics and Probability, Michigan State University, East Lansing, MI 48824}
\setcounter{Maxaffil}{0}
\renewcommand\Affilfont{\itshape\small}

\date {}

\maketitle

\begin{abstract}
Stochastic parameterizations are used in numerical weather prediction and climate modeling to help capture the uncertainty in the simulations and improve their statistical properties. Convergence issues can arise when time integration methods originally developed for deterministic differential equations are applied naively to stochastic problems. \cite{hodyss2013,hodyss2014} demonstrated that a correction term to various deterministic numerical schemes, known in stochastic analysis as the It\^o correction, can help improve solution accuracy and ensure convergence to the physically relevant solution without substantial computational overhead. The usual formulation of the It\^o correction is valid only when the stochasticity is represented by {\it white} noise. In this study, a generalized formulation of the It\^o correction is derived for noises of any color.  The formulation is applied to a test problem described by an advection-diffusion equation forced with a spectrum of fast processes. We present numerical results for cases with both constant and spatially varying advection velocities to show that, for the same time step sizes, the introduction of the generalized It\^o correction helps to substantially reduce time integration error and significantly improve the convergence rate of the numerical solutions when the forcing term in the governing equation is rough (fast varying); alternatively, for the same target accuracy, the generalized It\^o correction allows for the use of significantly longer time steps and hence helps to reduce the computational cost of the numerical simulation.

\end{abstract}

\section{Introduction}\label{introduction}

\textcolor{black}{Physical and chemical processes happening in the Earth's atmosphere span many orders of magnitude in terms of their spatial and temporal scales, which presents great challenges to numerical modeling. For example, in general circulation models, motions or phenomena that can not be resolved in space or time but have significant impact on the large-scale flow motions need to be accounted for using parameterizations (see e.g. \cite{mcfarlane2011}.)}

In recent years, stochastic parameterizations have become an active area of research (see review by \cite{berner2017,leutbecher2017}). The fundamental principle behind the stochastic formulation is that the state of the unresolved processes at any instant is not entirely determined by the state of the resolved processes. Thus, an element of randomness needs to be introduced to account for this indeterminacy. This randomness can act as a source of roughness in the temporal evolution of the governing equations' right-hand-side terms as well as in the evolution of the solution. Deterministic time integration
schemes used in numerical weather prediction and climate projection models, 
however, typically 
assume temporal smoothness of the underlying solutions. When such schemes are applied naively to stochastic parameterizations, the conditions for solution convergence might no longer be satisfied.

\textcolor{black}{
There is a large body of work dedicated to the development of stochastic numerical schemes and stochastic versions of deterministic numerical schemes especially targeted for weather and climate applications (see e.g. \cite{ruemelin1982, penland2002,sardeshmukh_2003,ewald2004,ewald2005,hansen2006}). The purpose of our work is more modest. We want to study the use of a correction term that can help improve the solution accuracy of deterministic schemes when part of the variables of the model is replaced by noise and in particular colored noise.}

As shown in \cite{hodyss2013}, the convergence issue of deterministic numerical schemes when applied to stochastic parameterizations can be investigated through the use of tools from stochastic analysis  \textcolor{black}{(see e.g Section 3.3 in \cite{oksendal2003}} \textcolor{black}{and Section 4.9 in \cite{kloeden1992}}). In particular, for the cases when an unresolved process is replaced by a rough random  process (e.g. white noise), it is not difficult to construct examples for which popular deterministic numerical schemes 
(e.g, Euler forward and backward, \textcolor{black}{Adams-Bashforth}) will no longer converge to the physically relevant solution except for special cases (e.g., the second-order Runge-Kutta scheme analyzed by \cite{hodyss2013}).
Multiple examples relevant for atmospheric modeling can be found in \textcolor{black}{Figure~3 in} the paper of \cite{hodyss2013}.
Here, ``physically relevant solution" refers to the one corresponding to ordinary calculus (see discussion below).
\textcolor{black}{In the study by \cite{hodyss2014}, ensemble simulations of Hurricane Isaac in the year 2012 were conducted using the Navy Operational Global Atmospheric Prediction System; it was shown that the choice of numerical scheme for the stochastic term can lead to failure in predicting the correct ensemble mean of  hurricane intensity.}

The mathematical reason for the lack of convergence is that when we replace an unresolved process with white noise, the equations describing the phenomena under investigation make sense only in integral form (not in the usual differential form). The integral form of the equations contains a temporal integral of an expression involving the white noise process. If we try to estimate such an integral through a limiting process involving progressively refined subintervals, different answers will be obtained depending on the manner we choose to discretize the interval of integration 
\textcolor{black}{(see e.g. Section 4.9 in \cite{kloeden1992})}. 
%
%
\textcolor{black}{The two most-studied discretization methods in stochastic analysis}
are: i) using the left endpoint of each subinterval (which leads to the It\^o integral or It\^o interpretation) and ii) using the middle point of each subinterval (which leads to the Stratonovich integral or Stratonovich interpretation). \textcolor{black}{The Stratonovich interpretation leads to ordinary calculus while the It\^o interpretation does not (see e.g. Sections 3.9 and 4 in \cite{oksendal2003}). Recall from the previous paragraph, that the physically relevant solution for the systems the weather and climate researchers are attempting to model is the one corresponding to ordinary calculus. Thus, the solution resulting from the Stratonovich interpretation is the physically relevant one \cite{hodyss2013}}.


It is important to note that many popular time integration schemes designed for deterministic problems will converge to the It\^o solution when applied to stochastic problems driven by {\it white} noise \cite{kloeden1992}. In other words, naively describing an unresolved process by  {\it white} noise and solving the stochastic equation with a deterministic numerical scheme can lead to erroneous results even in the limit of infinite temporal resolution. 
Fortunately, the It\^o and Stratonovich interpretations are related, and this relationship 
can help recover, at least to some extent, the convergence of deterministic numerical schemes to the physically relevant Stratonovich solution. The connection between the two interpretations comes in the form of a correction term called the It\^o correction. When the It\^o correction is added to the equation, the numerical solution under the It\^o interpretation converges to the Stratonovich solution.

While the above-mentioned It\^o-Stratonovich correspondence is a basic concept in stochastic analysis, the widely known form of the It\^o correction applies only to the case of white noise. A key feature of white noise is that it has zero auto-correlation (and hence no memory). Given the typical time step size of seconds to an hour in weather and climate models, some parameterized processes (e.g., turbulence and cumulus convection) can have characteristic time scales equivalent to multiple time steps. 
\textcolor{black}{Therefore,} colored noise, which has non-zero autocorrelation length, can provide a better description of such processes. 
\textcolor{black}{In fact, the state of the art in accounting for model uncertainties of Earth systems points to the need of stochastic processes with spatio-temporal correlations (see e.g. Sections 3 and 5 in \cite{leutbecher2017}), which makes our construction more relevant for applications.}  

A fundamental difference between colored noise and white noise is that colored noise is in principle resolvable while white noise is not. In other words, 
{\it if one could use small enough step sizes}, there would be no distinction between the It\^o and Stratonovich interpretations for the case of colored noise. All deterministic numerical schemes will eventually converge to the Stratonovich solution. But even in simple examples, let alone the very complex and expensive systems encountered in weather and climate prediction, the critical timestep that recovers convergence to the Stratonovich solution can be prohibitively small. As a result, for realistically affordable time step sizes, the dichotomy between It\^o and Stratonovich interpretations practically exists and needs to be addressed also for the cases of colored noise. \textcolor{black}{In other words, the long time steps in practical applications motivate us to find time-stepping methods with higher accuracy.} 

Another point worth mentioning is that, as \cite{hodyss2013} and \cite{hodyss2014} have pointed out, certain deterministic numerical schemes (e.g., the second-order Runge-Kutta scheme) have (an unaveraged version of) the It\^o correction already ``built-in" and hence perform better for stochastic problems. Such schemes are typically multi-stage schemes which require multiple evaluations of the right-hand side of the governing equations, making them very expensive for weather and climate models. The It\^o correction, in contrast, allows for the use of 
\textcolor{black}{single-stage schemes (e.g. the Euler forward scheme)}
to be complemented by a correction constructed only for the stochastic term, and hence can be cost-effective. 
\textcolor{black}{It should be noted that although the It\^o integral can be viewed as applying the Euler forward scheme to the stochastic term, other discretization methods can also fail to converge to the Stratonovich integral, see Section~3 of \cite{hodyss2013}. This means that the addition of an It\^o correction could aid in the restoration and/or acceleration of convergence to the Stratonovich solution also in the case of other numerical schemes.}


For these reasons, we present in this paper a generalization of the It\^o correction that is valid for noises of any color. We use an advection-diffusion equation with constant or spatially varying advection velocity to demonstrate that, for both white and colored noises,
the generalized It\^o correction can accelerate convergence to the Stratonovich solution when added to the Euler forward scheme. We demonstrate that 
improved convergence means higher accuracy for the 
same step size or, alternatively, larger step size 
(and hence lower computational cost) for the same target accuracy.
\textcolor{black}{These results from the simple but relevant test problem provide a proof of concept, which motivates further exploration of the generalized It\^o correction for the purpose of helping improve the solution accuracy and efficiency in atmospheric models of various complexity, including general circulation models using stochastic parameterizations.}

The remainder of the paper is organized as follows: Section \ref{derivation} presents the derivation of the generalized It\^o correction. Section \ref{test_problem} contains a presentation of the test problem, the advection-diffusion equation with constant and spatially varying advection velocity, along with analytical results (supplemental details can be found in the Appendix). Section \ref{results} contains numerical results. Finally, Section \ref{conclusions} contains a discussion of our results as well as suggestions 
for future work.



\section{The generalized It\^o correction}\label{derivation}




We consider the following deterministic differential equation 
\begin{equation}
\frac{\partial u}{\partial t} = D(u) + P(u),
\label{eq:pde_DP}
\end{equation}
where $D(u)$ and $P(u)$ are the resolved dynamics and 
parameterized physics, respectively. 
Here we focus on the special case where $P(u)$ takes the form 
\begin{equation}
P(u) = g(u) H(t).
\label{eq:P_deterministic}  
\end{equation}
This form results from the attempt to eliminate a fast-evoling physical quantity
from the original equations and replace it by a time-dependent process. We can consider a more general form where $H$ depends also on the spatial variable, but that generalization will not alter the derivation of the generalized It\^o correction below, hence we restrict our attention to the case where $H$ depends only on $t$.

If the time scales associated with ${H}(t)$ are substantially 
shorter than the time scales of ${D}(u)$, we can approximate 
${P}(u)$ by its stochastic counterpart ${P}_s(u)$ defined as
\begin{equation}
{P}_s(u) = g(u)\dot{{R}}(t),
\label{eq:P_stochastic}  
\end{equation}
where $\dot{{R}}$(t) represents a general noise term. \textcolor{black}{We note here that replacing $H(t)$ by the noise term, $\dot{{R}}(t),$ may include a limiting process where the function $g(u)$ may be also modified (see e.g. \cite{papanicolaou1974}).This does not alter the main line of our derivation and we keep the notation $g(u)$ for the multiplicative factor.}

\textcolor{black}{Using Eq.~\eqref{eq:P_stochastic}, we get the following stochastic counterpart of the deterministic equation originally given by Eq.~\eqref{eq:pde_DP}:}
\begin{equation}
\frac{\partial u}{\partial t} = {D}(u) + {P}_s(u).
\label{eq:sde_DP}
\end{equation}
Without loss of generality, we assume $\mathbb{E}[\dot{{R}}(t)] \equiv 0$, 
where $\mathbb{E}[\cdot]$ denotes
the mean over different realizations of the noise process. 
If $\mathbb{E}[\dot{{R}}(t)\dot{{R}}(t')] = \delta(t - t')$ 
where $\delta(\cdot)$ is Dirac's delta function, 
then $\dot{{R}}(t)$ is white noise and ${R}(t)$ is a Wiener process;
when $\mathbb{E}[\dot{{R}}(t)\dot{{R}}(t')] \neq \delta(t - t')$,
$\dot{{R}}(t)$ is a colored noise.

We focus on how to numerically
solve Eq.~\eqref{eq:sde_DP} after 
its form has been derived;
how to construct a good $P_s(u)$
to approximate the original $P(u)$ is a separate 
topic which is outside the scope of the current work.


%
%
%
%

\subsection{Derivation}
\label{sec:derive}

Let us take the integral over an arbitrary time window $(t_1, t_2)$ on both sides of Eq.~\eqref{eq:sde_DP}.
For ${P}_s(u)$, we discretize the time interval into $J$ bins of 
equal length $\Delta t$ and denote 
the increment of $R$ in each bin as $\Delta R$.
We use $t_j^*$ to denote the discretization point inside the $j$-th bin, 
i.e., the instant where the value of ${P}_s(u)$ is evaluated for numerical integration.
With this notation, the integral of Eq.~\eqref{eq:sde_DP} can be written as
\begin{eqnarray}
u(t_2) - u(t_1) 
&=&
\int_{t_1}^{t_2} D\big[u(t)\big] dt + 
\lim_{\Delta t\rightarrow 0}\sum_j g\big[u(t_j^*)\big]\Delta {R}_j\,. 
\label{eq:u_increment_generic}
\end{eqnarray}

In the white noise case (i.e., ${R}(t)=B(t)$ where $B(t)$ is the Wiener process), the choice of discretization point for the integral can lead to different results \cite{oksendal2003}. The two most popular choices are defined as
\begin{eqnarray}
\mbox{It\^o integral: } 
\int g(u) dB =
\lim_{\Delta t\rightarrow 0}\sum_j g\big[u(t_j^*)\big]\Delta B_j \notag \\
\text{where}\quad t_j^*=t_j \text{ (left endpoint),}
\label{eq:Ito_integral}
\end{eqnarray}
and 
\begin{eqnarray}
\mbox{Stratonovich integral: }
\int g(u)\circ dB =
\lim_{\Delta t\rightarrow 0}\sum_j g\big[u(t_j^*)\big]\Delta B_j 
\nonumber\\ 
\text{where}\quad  t_j^* =  \frac{t_j+ t_{j+1}}{2} = t_j + \frac{\Delta t}{2} \text{ (midpoint).} 
\label{eq:stratonovich_integral}
\end{eqnarray}
Here
$\Delta B_j = B_{j+1} - B_j .$ Because the physical processes represented by the deterministic equation \eqref{eq:pde_DP}  
are assumed continuous, the Stratonovich integral should be used in our case (see Section 3.3 in \cite{oksendal2003}). 

\textcolor{black}{Before proceeding further with the derivation, we note that the It\^o interpretation for the stochastic integral as shown in Eq. \eqref{eq:Ito_integral} coincides with how the Euler forward approach would treat the stochastic term, but the correspondence between It\^o and Stratonovich interpretations that we derive below is not tied to the forward Euler scheme.
We will come back to this point in Section~\ref{derivation}\ref{sec:other_time-stepping_schemes}.}

It is well known in stochastic analysis that in the white noise case, 
the Stratonovich integral can be written as the sum of an It\^o integral and a correction term called the It\^o correction (see, e.g., \cite{oksendal2003}). 
Below we show that the same is true for colored noise,
although the It\^o correction needs to be generalized.

For $t_j^* = \left(t_j+ t_{j+1}\right)/2,$ performing a Taylor expansion of 
$g\big[u(t_j^*)\big]$ about $t_j$ and 
expressing $\partial u/\partial t$ using Eq.~\eqref{eq:sde_DP} gives
\begin{eqnarray}
g\big[u(t_j^*)\big] &=& 
g\big[u(t_j)\big] + 
\frac{\Delta t}{2}
\left( \frac{dg(u)}{du} \frac{\partial u}{\partial t} \right) \!\!\bigg\vert_{t_j}  
+ \frac{\Delta t^2}{8} \left(\frac{d^2g[u(t)]}{dt^2}\right)\!\!\bigg\vert_{\xi}
\\
&=& g\big[u(t_j)\big] + 
  \left(\frac 12 \frac{dg(u)}{du} {D}[u] \right)\!\!\bigg\vert_{t_j}\!\!\Delta t
+ \left(\frac 12 \frac{dg(u)}{du} g[u]\dot{{R}}(t)\right) \!\!\bigg\vert_{t_j} \!\!\Delta t
\nonumber\\
&& + \frac{\Delta t^2}{8} \left(\frac{d^2g[u(t)]}{dt^2}\right)\!\!\bigg\vert_{\xi}
\label{eq:g_t*_taylor}
\end{eqnarray}
where $\xi \in \left[t_j, \left(t_j+ t_{j+1}\right)/2\right]$.
For small $\Delta t$, we write
\begin{equation}
    \dot{{R}}(t_j) \Delta t \approx \Delta {R}_j. 
\end{equation}
%
Hence, Eq.~\eqref{eq:g_t*_taylor} can be approximated as 
\begin{gather}
g\big[u(t_j^*)\big] \approx g\big[u(t_j)\big]  
+ \left(\frac 12 \frac{dg(u)}{du} {D}(u) \right)\!\!\bigg\vert_{t_j}\!\!\Delta t
+ \left(\frac 12 \frac{dg(u)}{du} g[u]\right) \!\!\bigg\vert_{t_j} \!\!\Delta {R}_j \notag \\
+ \frac{\Delta t^2}{8} \left(\frac{d^2g[u(t)]}{dt^2}\right)\!\!\bigg\vert_{\xi} \label{eq:g_t*_taylor_cont}
\end{gather}

Assuming $g(u)$ is sufficiently smooth and $\Delta t$ is small, 
one can neglect the second and fourth terms 
on the right-hand side of Eq.~\eqref{eq:g_t*_taylor_cont} but, in general, not the third term.  Therefore, \textcolor{black}{with the Stratonovich interpretation of the stochastic integral in Eq.~\eqref{eq:u_increment_generic}, we have}

\begin{gather}
	u(t_2) - u(t_1) =
	    \underbrace{
	    \int_{t_1}^{t_2} {D}(u)dt
	    }_{\mbox{traditional integral}}   
	   + \underbrace{
	    \int_{t_1}^{t_2} g(u)d{R}
	    }_{\mbox{It\^o integral}}  \notag \\
	   + \lim_{\Delta t\rightarrow 0}\sum_j 
        \left(\frac 12 \frac{dg(u)}{du} g[u]\right) \!\!\bigg\vert_{t_j}\!\!
		\left(\Delta{R}_j\right)^2 \label{eq:u_increment_stratonovich}
\end{gather}
The mathematical expectation of the last term in Eq.~\eqref{eq:u_increment_stratonovich} is
\begin{eqnarray}
		\lim_{\Delta t\rightarrow 0}\sum_j 
        \left(\frac 12 \frac{dg(u)}{du} g[u]\right) \!\!\bigg\vert_{t_j}\!\!
		\mathbb{E}\big[ \left(\Delta{R}_j\right)^2 \big]\,,
		\label{eq:Ito_generalized}
\end{eqnarray}
which is the generalized It\^o correction in its integral form.
The exact form of the expectation in expression~\eqref{eq:Ito_generalized} depends on the formulation of ${R}.$
For example, \textcolor{black}{when $R$ is the Wiener process, the increment $\Delta R_j$ is a Gaussian random variable with mean 0 and variance $\Delta t$, i.e.,}
\begin{equation}
	\mathbb{E}\big[ \left(\Delta{R}_j\right)^2 \big] = \Delta t ,
\end{equation}
hence \eqref{eq:Ito_generalized} becomes 
\begin{equation}
	    \int_{t_1}^{t_2} \left(\frac 12 \frac{dg(u)}{du} g(u)\right) dt ,
\end{equation}
which is the integral form of the traditional It\^o correction (see, e.g. Section 3.3 in \cite{oksendal2003}).

\textcolor{black}{The generalized It\^o correction \eqref{eq:Ito_generalized} can be extended to the case of multiple partial differential equations (PDEs)  each containing multiple noise processes. Let us assume a system of $n$ PDEs for the functions ${\bf u}=(u_1,u_2,\ldots,u_n)$:
%
%
\begin{equation}
   \frac{\partial u_i}{\partial t} = D_i(u) + \sum_{l=1}^p g_{il}(u)\dot{R}_l(t), \; \text{for} \; i=1,\ldots,n
\end{equation}
%
where $\dot{R_l}(t)=(\dot{R_1}(t),\dot{R_2}(t),\ldots,\dot{R_p}(t))$ is a $p$-dimensional vector noise process with independent components.
%
Then, the expression for the generalized It\^o correction for the equation for $u_i$ is given by
\begin{eqnarray}
		\lim_{\Delta t\rightarrow 0}\sum_j 
        \left(\frac 12  \sum_{l=1}^p \sum_{k=1}^n \frac{\partial g_{il}(u)}{\partial u_k} g_{kl}[u]\right) \!\!\bigg\vert_{t_j}\!\!
		\mathbb{E}\big[ \left(\Delta R_{lj}\right)^2 \big]\,.
		\label{eq:Ito_generalized_multiple}
\end{eqnarray}
where $\Delta R_{lj} = R_{l,j+1}- R_{lj}$ is the increment of the $l$th noise process $R_l.$
}

\subsection{Remarks}
\label{sec:remarks}

We want to make \textcolor{black}{two} remarks concerning the derivation of the generalized It\^o correction \eqref{eq:Ito_generalized}. First, there is an alternative way to derive the generalized It\^o correction. In particular, under the assumption that the correlation time of the noise is short, one can employ the expansion devised by Stratonovich (see Section 4.8 in \cite{stratonovich1963}), through which a stochastic equation driven by colored noise can be rewritten as an {\it effective} stochastic equation driven by white noise. Then, one can compute the traditional It\^o correction for the resulting white noise driven equation. 

Second, the It\^o correction, in its traditional or generalized form, can be interpreted as a memory term encountered in model reduction formalisms (see e.g. \cite{chorin2006problem}). By using as discretization point the left endpoint of each interval, the It\^o interpretation of the stochastic integral makes the evolution of the stochastic process $\dot{{R}}(t)$ {\it independent} of the solution $u(t).$ The It\^o correction serves as a way to account for the interaction of $\dot{{R}}(t)$ and $u(t)$ during the interval $\Delta t,$ similar to the role played by memory terms in model reduction which account for the interaction between resolved and unresolved variables.  

\subsection{\textcolor{black}{Applicability}}
\label{sec:other_time-stepping_schemes}

\begin{color}{black}
It has been stated earlier in Section \ref{derivation}\ref{sec:derive} 
that our generalized It\^o correction \eqref{eq:Ito_generalized} 
is not tied to the specific discretization method (e.g., Euler forward) 
that is chosen 
for the time integral of the stochastic term in Eq.~\eqref{eq:sde_DP}.
The reason is that for any discretization, 
as long as an analysis similar to 
Eqs.~\eqref{eq:g_t*_taylor}--\eqref{eq:g_t*_taylor_cont} reveals that the discretized integral converges to 
the It\^o integral,
expression~\eqref{eq:Ito_generalized}  can be used to 
obtain numerical results that converge to the Stratonovich solution.

The method of analysis demonstrated by 
Eqs.~\eqref{eq:g_t*_taylor}--\eqref{eq:g_t*_taylor_cont} 
can also be applied to the ``decentered" time-stepping methods 
commonly used in atmospheric models.  Since these methods  
approximate time integrals (or derivatives) using
the discretization point
$t^*_j = (1-\lambda)t_j + \lambda t_{j+1} = t_j + \lambda\Delta t$ 
where $0 \le \lambda \le 1$, 
Eq.~\eqref{eq:u_increment_generic} becomes 
\begin{gather}
	u(t_2) - u(t_1) =
	    \underbrace{
	    \int_{t_1}^{t_2} {D}(u)dt
	    }_{\mbox{traditional integral}}   
	   + \underbrace{
	    \int_{t_1}^{t_2} g(u)d{R}
	    }_{\mbox{It\^o integral}} \notag \\
	   + \lim_{\Delta t\rightarrow 0}\sum_j 
        \left(\lambda \frac{dg(u)}{du} g[u]\right) \!\!\bigg\vert_{t_j}\!\!
		\left(\Delta{R}_j\right)^2 \,. \label{eq:u_increment_decentered}
\end{gather}
A comparison of Eq.~\eqref{eq:u_increment_decentered} with
Eq.~\eqref{eq:u_increment_stratonovich} suggests that 
the correction term linking the ``decentered" integral 
and the Stratonovich integral is
\begin{eqnarray}
		\lim_{\Delta t\rightarrow 0}\sum_j 
        \left((\frac 12-\lambda) \frac{dg(u)}{du} g[u]\right) \!\!\bigg\vert_{t_j}\!\!
		\mathbb{E}\big[ \left(\Delta{R}_j\right)^2 \big],
		\label{eq:Ito_generalized_decenter}
\end{eqnarray}
which is a further generalization of
expression~\eqref{eq:Ito_generalized} for non-zero $\lambda$.
In the special case of $\lambda=0$, the ``decentered" scheme 
gives the It\^o interpretation of the stochastic integral and we recover \eqref{eq:Ito_generalized}; the case of $\lambda=1/2$ corresponds to the Stratonovich interpretation of the integral and the correction vanishes. 
For the case when $R$ is the Wiener process, expression \eqref{eq:Ito_generalized_decenter} corresponds to the correction formula for {\it white} noise that is found in \cite{hodyss2013} (see equations~(2.7) and (3.21) therein) as well as in Section 3.5 of \cite{kloeden1992}. 

The example of decentered schemes discussed above can be generalized even further: 
for a generic discretization method, an analysis 
similar to Eqs.~\eqref{eq:g_t*_taylor}--\eqref{eq:u_increment_stratonovich}, 
followed by a comparison with the desired interpretation of the equation 
(i.e., Eq.~\eqref{eq:u_increment_stratonovich} for Stratonovich or 
Eq.~\eqref{eq:u_increment_stratonovich} without 
the last right-hand-side term for It\^o),
can lead to the correction term needed to obtain numerical results
converging to the desired type of solution (Stratonovich or It\^o).
\end{color}


\subsection{\textcolor{black}{Integral versus differential form}}
\label{sec:integral_vs_differential}
\textcolor{black}{
The expressions for the generalized It\^o correction 
presented so far have been obtained from the {\it integral} form of 
the stochastic equation, while in the literature on stochastic analysis,
the It\^o correction conventionally denoted by the symbol $I$ 
is typically the term that is added to the {\it differential} 
form of the stochastic equation. 
The differential form of our generalized It\^o correction 
in the test problem discussed below is given in 
Section~\ref{test_problem} and in Appendix~\ref{sec:noise}.
An example showing how the differential form
is obtained from the integral form can be found in 
Appendix~\ref{sec:noise} (Eqs.~\ref{eq:colored_noise_variance_2}--\ref{ito_correction_1}).
}



\section{Test problem}\label{test_problem}
  In the remainder of the paper, we use an example to 
demonstrate the impact of the generalized It\^o correction.
We consider the following stochastic differential equation
\begin{equation}\label{advection_diffusion}
\frac{\partial u}{\partial t} = -\left[c+\frac{\epsilon}{2}\cos(x)\right]\frac{\partial u}{\partial x} + \mu \frac{\partial^2 u}{\partial x^2}+ g(u)n(t),
\end{equation} 
with initial condition $u(x,0) = u_0(x)$ and periodic boundary conditions on $[0,2\pi].$
In the context of atmosphere modeling, 
the first two terms on the right-hand side represent the resolved dynamics 
and the last term represents fast varying physics parameterizations.
When the parameter $\epsilon$ is set to 0, we recover the advection-diffusion equation 
with constant advection velocity discussed in \cite{hodyss2013}.
The inclusion of $\epsilon\cos(x)/2$ 
in the first right-hand-side term
makes the advection velocity spatially varying.
In Section~\ref{results}, numerical results are shown for
both $\epsilon=0$ and $\epsilon=10^{-3}$.
Following \cite{hodyss2013}, we let 
\begin{eqnarray}
&&c=1,\quad \mu=0.1,\, \\
&&g(u) = \rho \frac{\partial u}{\partial x\,} \,\,\,\text{with}\, \,\rho = 0.2\,.
\end{eqnarray}
The stochastic noise process $n(t)$ 
is the same as described in Appendix A of \cite{hodyss2013} 
(also described in Appendix~A of this paper).
For this choice of $g(u)$ and $n(t)$, the generalized It\^o correction is given by
\begin{eqnarray}
I = \frac 12 \rho^2 \frac{\partial^2 u}{\partial x^2} 
\lim_{N_f \rightarrow \infty}\frac{1}{N_f}\! \left[
 \frac{C(\omega_0)^2}{2} 
+ \sum_{m=1}^{N_f} C(\omega_m)^2
\right], \label{ito_correction}
\end{eqnarray}
where $N_f$, $C$, $\omega_0$ and $\omega_m$ are parameters of the 
noise process $n(t)$ (cf. Appendix~A). 
\textcolor{black}{The expression for $C$ contains a parameter 
$\alpha$ that controls the color of the Fourier spectrum of $n(t)$, with
$\alpha = 0$ corresponding to white noise, and larger $\alpha$ values
corresponding to noise spectra that are more red.}

\textcolor{black}{Before we proceed, we want to make an important remark about the formula for $I.$ For the case of white noise ($\alpha=0$), we obtain $$\lim_{N_f \rightarrow \infty}\frac{1}{N_f}\! \left[
 \frac{C(\omega_0)^2}{2} 
+ \sum_{m=1}^{N_f} C(\omega_m)^2
\right]=1,$$ and we recover the usual It\^o correction expression. However, for the case of colored noise  ($\alpha \neq 0$) with exponentially decaying spectrum, we have $\lim_{N_f \rightarrow \infty}\frac{1}{N_f}\! \left[
 \frac{C(\omega_0)^2}{2} 
+ \sum_{m=1}^{N_f} C(\omega_m)^2
\right]=0.$ This result is not surprising. As we have explained also in the introduction, in the case of colored noise, the distinction between the It\^o and Stratonovich interpretations disappears in the limit of infinite temporal resolution. The reason is that in the limit of infinite temporal resolution a non-white colored noise is resolved and thus all discretizations of the integral of the stochastic term give the same answer. However, {\it any} numerical experiment that one conducts always has {\it finite} temporal resolution. In this case, the generalized It\^o correction for the case of colored noise is no longer zero. Moreover, as we show with our numerical results, it can play a significant role in restoring or accelerating convergence to the Stratonovich solution. This remark is particularly pertinent for weather and climate applications where due to computational limitations we are always forced to use larger timesteps than the shortest timescales present in the solution.
}

To derive analytical solutions for the test problem, 
we express $u(x,t)$ in the form of a superposition of 
Fourier modes
\begin{equation}\label{eq:solution_Fourier}
u(x,t) = \sum_{k\in \mathbb{Z}}F_k(t)\exp(ikx)
\end{equation}
and transform Eq.~\eqref{advection_diffusion} into a system of 
stochastic differential equations.
Here $i$ is the imaginary unit.
Like in \cite{hodyss2013}, we assume the initial condition 
contains only one mode, i.e.,
\begin{equation}\label{eq:initial_condition}
u(x,0) = \cos(k_0x)\quad\text{with}\quad k_0 = 1.
\end{equation}

\subsection{Case with constant advection velocity ($\epsilon=0$)}

As pointed out by \cite{hodyss2013}, when the advection velocity is constant, the Fourier modes
are uncoupled. 
The ordinary differential equation (ODE) for $F_k(t)$ reads
\begin{equation}
\frac{d F_k}{d t} = -ik c F_k - \mu k^2 F_k + ik \rho n(t)F_k \,.
\label{eq:fourier_adv_diff_constant}
\end{equation}
The analytical solution of Eq.~\eqref{eq:fourier_adv_diff_constant} takes the form
\begin{equation}\label{eq:analytic_solution_constant_coeff}
F_k(t) = A\exp\left[-\left(ick+\mu k^2\right)t + i\rho k\int_0^t n(t')dt'\right]
\end{equation}
with $A$ being any complex constant. 
Initial condition \eqref{eq:initial_condition} implies that 
$A = 1$ in Eq.~\eqref{eq:analytic_solution_constant_coeff};
only one Fourier mode (the one corresponding to $k=1$) is sufficient 
to represent the solution (the Fourier mode for $k=-1$ is also needed but due to the solution of \eqref{advection_diffusion} being real, it is the complex conjugate of the solution for the Fourier mode with $k=1.$)

\subsection{Case with spatially varying advection velocity ($\epsilon \neq 0$)}

For cases with nonzero $\epsilon$, even if the initial condition has a single Fourier mode, 
the spatially dependent component of the advection velocity causes the representation 
of the solution to require more than one mode. 
This is an elementary way to introduce coupling between different Fourier modes
but still keep Eq.~\eqref{advection_diffusion} linear.

We truncate the Fourier series in Eq.~\eqref{eq:solution_Fourier}
to retain only modes with appreciable magnitudes 
and denote the largest remaining wavenumber as $N_x$. Substituting Eq.~\eqref{eq:solution_Fourier} for the unknown $u$ 
in \eqref{advection_diffusion} gives
\begin{equation}\label{eq:Fourier_ODE_nonconstant}
\begin{split}
\sum_{k=-N_x}^{N_x}\frac{d F_k(t)}{d t}\exp(ikx) =&
\sum_{k=-N_x}^{N_x}\left\{-ik{F_k(t)}
\exp(ikx)\left[c+\frac{\epsilon}{2}\cos(x)\right]\right. \\
&\left.- \mu k^2{F_k(t)}\exp(ikx)+ ik\rho F_k(t)\exp(ikx)n(t) \frac{}{}\right\} \,.
\end{split}
\end{equation} 
By multiplying Eq.~\eqref{eq:Fourier_ODE_nonconstant} 
with $\exp(-ikx)$ and integrating over $[0,2\pi]$, 
we get the following coupled equations:
\begin{itemize}
\item for $k = -N_x+1,\dots, N_x-1$,
\begin{gather}\label{eq:n}
\frac{d F_k(t)}{d t} = (-ick-\mu k^2)F_k(t)-\frac{i(k+1)\epsilon}{4}F_{k+1}(t)-\frac{i(k-1)\epsilon}{4}F_{k-1}(t) \\
+ik\rho F_k(t)n(t)\,;
\end{gather}
\item for $k = N_x$, 
\begin{equation}\label{eq:N_x}
\frac{d F_k(t)}{d t} =(-ick-\mu k^2)F_k(t)-\frac{i(k-1)\epsilon}{4}F_{k-1}(t) +ik\rho F_k(t)n(t); 
\end{equation}
\item for $k = -N_x$,
\begin{equation}\label{eq:-N_x}
\frac{d F_k(t)}{d t} = 
(-ick-\mu k^2)F_k(t)-\frac{i(k+1)\epsilon}{4}F_{k+1}(t) +ik\rho F_k(t)n(t)\,. 
\end{equation}

\end{itemize}
Using the notation defined in Appendix~B, we can write the above stochastic ODE system 
for the Fourier mode coefficients $F_k$ in matrix form as
\begin{equation} \label{eq:F_equation_matrix_form}
\frac{d \mathbf{F}}{d t} = \left[ \pmb{\mathsf{D}} + \rho n(t) \pmb{\mathsf{H}} \right]\mathbf{F} \,.
\end{equation} 
The analytical solution reads 
\begin{equation}\label{eq:analytic_solution_nonconstant}
\mathbf{F}(t) =  \exp\left(\pmb{\mathsf{D}}\,t+\pmb{\mathsf{H}}\rho\int_0^t n(t')dt' \right)\mathbf{F}(0) .
\end{equation}

\section{Numerical Results}\label{results} 








In this section we use numerical results to show how noise $n(t)$ of different color (roughness in time) can affect
\textcolor{black}{the solution convergence.}
\textcolor{black}{We also} demonstrate how the inclusion of the generalized It\^o correction can help in restoring and/or accelerating convergence. 
\textcolor{black}{As explained in Section~\ref{derivation}\ref{sec:other_time-stepping_schemes}, 
the validity of our generalized It\^o correction is 
not tied to any specific time-stepping method. 
For simplicity, we use in our numerical experiments the forward Euler scheme as an illustrating example.}

\subsection{Definition of solution error}

The error of a numerical solution is evaluated after two time units 
of integration using the $L_2$ norm (\cite{hodyss2013} and personal communication):
\begin{equation}\label{eqn:solution_error}
\mathcal{E}(\Delta t) = 
\left\{
\int_{0}^{2\pi} \left[ \widehat{u}(x,t=2)-u(x,t=2) \right]^2 dx
\right\}^{\frac 12} \,.
\end{equation}
Here $\widehat{u}$ and $u$ are the discrete and analytical solutions, respectively.
To ensure the accuracy of the analytical solution computed for our error evaluation,
the time integral of the noise process in 
Eqs.~\eqref{eq:analytic_solution_constant_coeff}
and \eqref{eq:analytic_solution_nonconstant} 
is calculated analytically.

\subsection{Case with constant advection velocity ($\epsilon=0$)}

For the case with $\epsilon$~=~0,
the discretization of Eq. \eqref{eq:fourier_adv_diff_constant} using forward Euler 
with the It\^o correction included is given by 
\begin{equation}
\dfrac{\widehat{F}_k(t_{j+1}) - \widehat{F}_k(t_{j})}{\Delta t}
  =  - ik c \widehat{F}_k(t_j)
     - \mu k^2 \widehat{F}_k(t_j) 
     + ik \rho \widehat{F}_k(t_j) n(t_j)  + I_k(t_j),
\end{equation}
where $n(t_j)$ is the colored noise at $t = t_j$ and $I_k(t_j)$ 
is the It\^o correction for $F_k$ at $t_j$,
\begin{equation}
I_k(t_j) = -\frac{1}{2} \rho^2 k^2 \widehat{F}_k(t_j) \frac{1}{N_f}\left\{
 \frac{C(\omega_0)^2}{2} 
+ \sum_{m=1}^{N_f} C(\omega_m)^2
\right\}. \label{ito_correction_discrete}   
\end{equation}
\textcolor{black}{We note that although the expression for the It\^o correction in Eq. \eqref{ito_correction} involves a limiting process, the limit is not present in the expression in Eq. \eqref{ito_correction_discrete} because we have discretized the equation and thus have picked a finite timestep.}

\begin{figure}[ht]
   \centering
   \includegraphics[width = 0.9\textwidth]{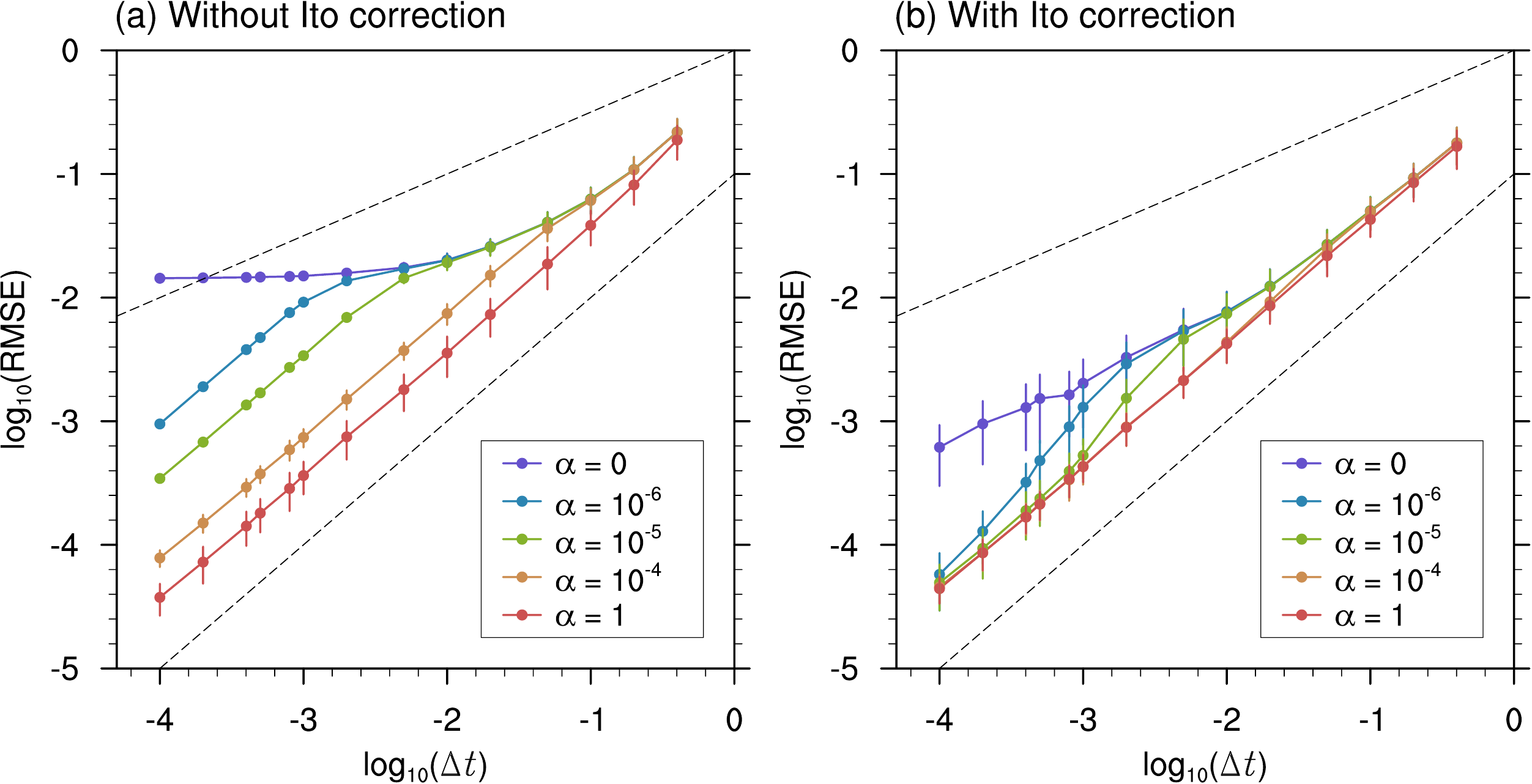}
   \caption{Error in the numerical solution of the 1D advection-diffusion equation with {\it constant} advection velocity ($\epsilon$~=~0 in Eq.~\ref{advection_diffusion}) and the dependency on time step size (x-axis) and characteristics of the noise term ($\alpha$~=~0, 10$^{-6}$, 10$^{-5}$, 10$^{-4}$, or 1, shown in different colors). The left and right panels show results obtained using the forward Euler scheme without and with the generalized It\^o correction, respectively. Simulations were performed for 100 realizations of the noise process and the $l_2$ solution error was calculated separately for each realization using Eq.~\eqref{eqn:solution_error}. The thick dots are the mean error of the 100 realizations; the vertical bars denote the standard deviation around the mean. \textcolor{black}{The two dashed lines are reference lines indicating convergence rates of 0.5 (upper) and 1.0 (lower), respectively.} }
   \label{figure_constant}
\end{figure}

Panel (a) of Figure \ref{figure_constant} (which appears also in \cite{hodyss2013}) shows the effect of different noises on the convergence of the Euler scheme {\it without} the It\^o correction. The thick dots are the $l_2$ error of the numerical solution averaged over 100 realizations of the noise process; the error bars denote the standard deviation around the average. Using the terminology of stochastic analysis, this plot (and the rest of them in the paper) shows the {\it strong} convergence of the numerical solution.
\footnote{As a reminder, we note that strong convergence is measured by the mean of the solution error of individual realizations of the stochastic equation while weak convergence is measured by the error of the mean solution.
}

We make two observations. First, for the case of white noise ($\alpha=0$, purple line), the Euler scheme {\it without} the It\^o correction fails to converge to the analytical solution no matter how small the step size is (it converges to the It\^o solution, cf. \textcolor{black}{Section 10.2 in}~\cite{kloeden1992}). Second, for the case of colored noise ($\alpha \neq 0$, blue, green, orange and red lines), the Euler scheme {\it without} It\^o correction will start converging to the analytical solution with order 1 (as predicted by deterministic numerical analysis, \textcolor{black}{see e.g. Chapter I.7} in \cite{hairer1993}) when the step size becomes smaller than some critical step size which depends on the color of the noise (value of $\alpha$). The more red the noise is (larger $\alpha$), the larger is the critical step size (see also \cite{hodyss2013} for a discussion and estimation of the critical stepsize).    

The right panel in Figure \ref{figure_constant} shows the effect of including the It\^o correction. We want to make again two observations. First, for the case of white noise ($\alpha=0$, purple line), the Euler scheme {\it with} the It\^o correction {\it does} converge to the analytical solution with order $1/2$ (see \cite{hu2016} for an explanation of this convergence rate). We note that this numerical result was mentioned in \cite{hodyss2013} although not illustrated by any graphic there. Second, for the case of colored noise ($\alpha \neq 0$, blue, green, orange and red lines), the Euler scheme with the {\it generalized} It\^o correction starts converging to the analytical solution with order 1 for larger step sizes than the Euler scheme {\it without} the It\^o correction. Thus, the addition of the It\^o correction can help restore and/or accelerate convergence of the forward Euler scheme.

\subsection{Case with spatially varying advection velocity ($\epsilon \neq 0$)}

We continue with the case of a spatially-dependent advection velocity with $\epsilon=10^{-3}.$ A small value was chosen for $\epsilon$ because the forward Euler scheme is explicit and only first-order. As such, it needs a very large number of steps in order to reach the asymptotic convergence regime for larger values of $\epsilon$ due to the need to resolve steepening gradients associated with the oscillatory nature of the spatial perturbation of the advection velocity. Moreover, the cost of evaluation of the noise $n(t),$ which depends quadratically on the number of timesteps, becomes very large when $\epsilon$ is large. For practical purposes (computational cost), we chose a small $\epsilon$  for the demonstration here.

\begin{figure}[ht]
   \centering
   \includegraphics[width = 0.9\textwidth]{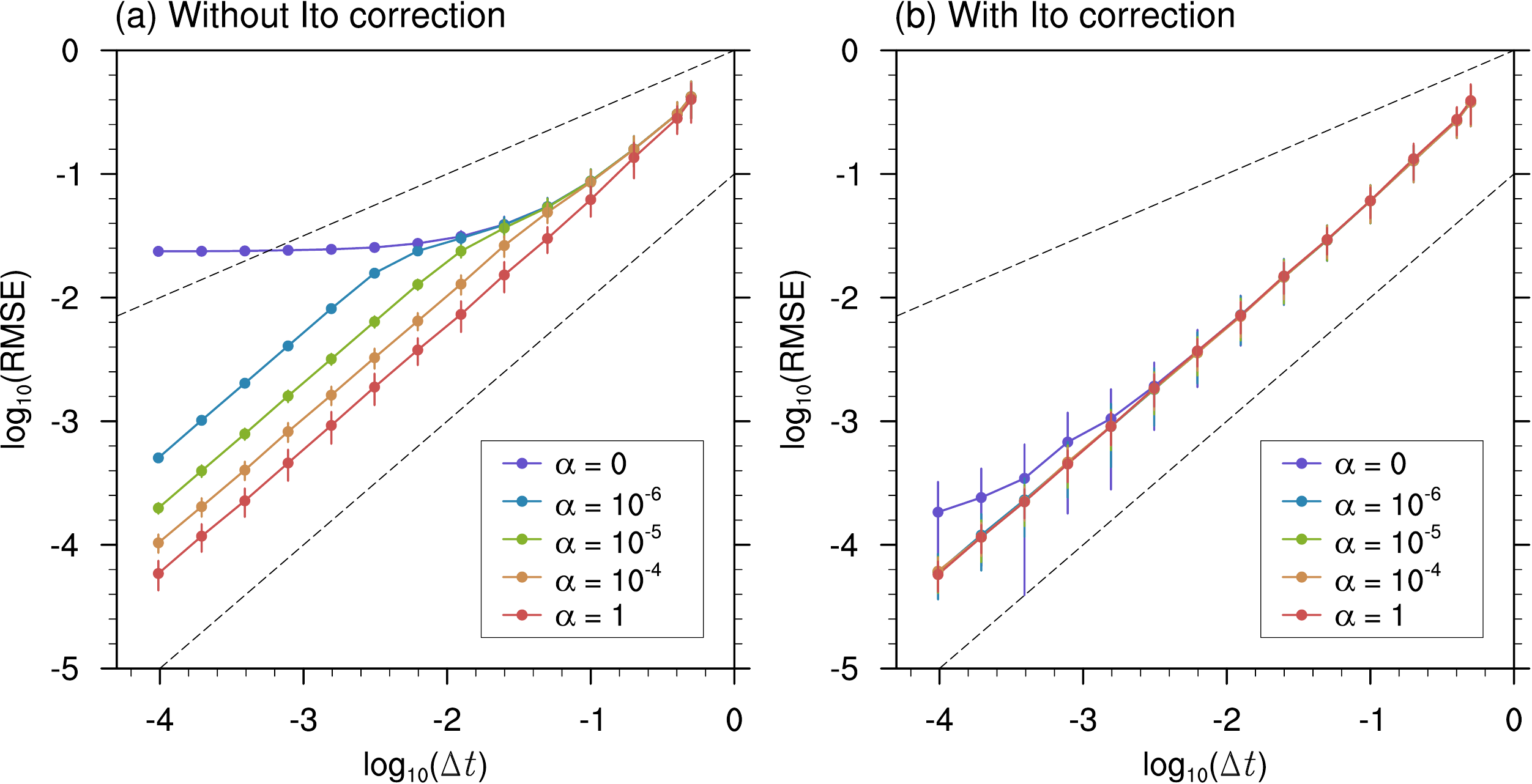}
   \caption{As in Figure~\ref{figure_constant} but for the case of $\epsilon$~=~10$^{-3}$ in Eq.~\eqref{advection_diffusion}.}
   \label{figure_dependent}
\end{figure}

We discretized Eq.~\eqref{eq:F_equation_matrix_form} using the Euler scheme with the generalized It\^o correction, i.e., 
\begin{equation}\label{eq:Euler_forward_with_Ito_nonconstant_case}
\dfrac{\widehat{\mathbf{F}}(t_{j+1}) - \widehat{\mathbf{F}}(t_{j})}{\Delta t} =
 \pmb{\mathsf{D}} \widehat{\mathbf{F}}(t_{j}) 
+ \rho n(t_{j}) \pmb{\mathsf{H}}\widehat{\mathbf{F}}(t_{j})
+ \mathbf{I}(t_j)
\end{equation}
where the It\^o correction reads
\begin{equation}
\mathbf{I}(t_j) = \frac{\rho^2}{2}\frac{1}{N_f}\left(\frac{C(\omega_0)^2}{2}
+\sum_{m=1}^{N_f} C(\omega_m)^2\right) \pmb{\mathsf{G}} \widehat{\mathbf{F}}(t_j)
\end{equation}
with the matrix $\pmb{\mathsf{G}}$ being
\begin{equation}
\pmb{\mathsf{G}} = \textrm{Diag}\{-(-N_x)^2, -(-N_x+1)^2,\cdots,-(N_x-1)^2 , -N_x^2\}\,.
\end{equation} 
The truncation wavenumber $N_x$ was chosen empirically:
a test simulation was conducted using Eq.~\eqref{eq:Euler_forward_with_Ito_nonconstant_case}
with a large $N_x$; an inspection of the magnitude of the resulting $\widehat{F}_k$
revealed $N_x = 5$ was sufficient to retain all modes 
with $|\widehat{F}_k|> 10^{-4}$. Hence, $N_x = 5$ was then used 
to obtain the results shown in Figure~\ref{figure_dependent}.

Figure~\ref{figure_dependent} shows that for the case of white noise ($\alpha=0$, purple line), the forward Euler scheme {\it without} the It\^o correction fails to converge to the analytical solution as expected. Figure \ref{figure_dependent} shows how the inclusion of the It\^o correction can restore convergence with order $1/2.$ We note that the standard deviation bars around the mean appear larger than in the case with $\epsilon=0$ because of the logarithmic scale of the plot. This figure demonstrates that for the case of colored noise ($\alpha \neq 0$), the use of the generalized It\^o correction again accelerates the establishment of the order 1 convergence regime predicted by deterministic numerical analysis.



\section{Conclusions}\label{conclusions}





Stochastic parameterizations are increasing in popularity in numerical weather prediction and climate modeling as a way to improve the statistical representation of the studied phenomena. Naive implementation of such parameterizations with deterministic numerical time integration schemes can cause serious convergence issues. Such issues can be alleviated by the addition of certain correction terms (known as the It\^o correction in stochastic analysis) to the deterministic numerical schemes. However, the It\^o correction was originally derived only for the special case when the stochastic process is represented by white noise. For numerical weather prediction and climate modeling it will be useful to have the option to properly handle colored noise.

We have derived a generalized It\^o correction for the case of colored noise and applied it to a test problem of an advection-diffusion equation driven by noise of different colors. Our results indicate that the generalized Ito correction can  substantially reduce the time discretization error, significantly improve the convergence rate of the numerical solutions and allow for the use of significantly larger step sizes. 

While our derivation started from a stochastic differential equation, the fact that colored noise is in principle resolvable by sufficiently small step sizes implies that the generalized It\^o correction can also be useful for deterministic problems for the purpose of improving solution convergence, accuracy, and efficiency. Compared to higher-order schemes like the Runge-Kutta family, the It\^o correction is less costly in terms of computing time; compared to implicit methods that may provide better stability, the It\^o correction is less intrusive in terms of the code modification it requires.

In the future, we plan to apply the current framework to more realistic atmospheric modeling problems, e.g. simplified versions of the atmospheric general circulation models or their parameterizations. \textcolor{black}{We acknowledge the fact that the parameterizations will likely not be given directly in the multiplicative form required by our formulation. However, there is hope that new approaches e.g. training a neural network to represent the function $g(u)$ in Eq.~\eqref{eq:P_deterministic} will allow us to still use our construction for more involved parameterizations.}

\textcolor{black}{The generalized It\^o correction for the case of colored noise still assumes that the subgrid phenomena evolve on a significantly shorter timescale than the phenomena we resolve explicitly. However, as was mentioned in Section \ref{sec:remarks}, the It\^o correction can be interpreted as a memory term in model reduction formalisms. This opens the possibility of exploiting more sophisticated types of memory terms which correspond to more nuanced and realistic noise processes. 
}


  
\section*{Acknowledgements}
   The authors thank Drs. Christopher J. Vogl (LLNL), Carol S. Woodward (LLNL), and Shixuan Zhang (PNNL) for helpful discussions on the numerical examples shown in this paper, and Dr. Daniel Hodyss (NRL) for clarifications regarding his earlier work that inspired our study.
This work was supported by the U.S. Department of Energy (DOE), 
Office of Science, 
Office of Advanced Scientific Computing Research (ASCR) and 
Office of Biological and Environmental Research (BER),
Scientific Discovery through Advanced Computing (SciDAC) program.
Pacific Northwest National Laboratory is operated by 
Battelle Memorial Institute for DOE under Contract DE-AC05-76RL01830.

\appendix

\section*{Appendix A}
  
\label{sec:noise}

Following \cite{hodyss2013}, we use the following specification of the noise process $n(t)$:
\begin{gather}
n(t) = \frac{1}{\sqrt{N_f \Delta t}}
 \left(C(\omega_0)\frac{b_0}{\sqrt{2}} + \sum_{m=1}^{N_f} C(\omega_m) \big[a_m\sin(\omega_m t) + b_m\cos(\omega_m t)\big]\right) \label{noise_1} \\
C(\omega) =  e^{-\alpha \omega^2} \\
\omega_m = \dfrac{2\pi m}{(N-1)\Delta t} \label{eq:omega}\\
N_f = (N-1)/2
\end{gather}
where $N$ is the number of discrete time levels per unit time, 
including the starting and ending time levels. 
$\alpha$ is a parameter controlling the color of the Fourier spectrum of $n(t)$ \textcolor{black}{($\alpha=0$ corresponds to {\it white} noise while $\alpha \neq 0$ to {\it colored} noise)}. 
To construct different realizations of the noise process, 
we \textcolor{black}{sample, for $m=0,\ldots,N_f,$ the coefficients $a_m$ and $b_m$ independently} from the normal distribution $\mathcal{N}(0,1)$.  
It should be noted that $n(t)$ is an {\it approximate} random noise. The difference between $n(t)$ 
and the noise term  $\dot{{R}}(t)$ in 
Section~\ref{derivation} is that 
$\dot{{R}}(t)$ would contain an infinite number of Fourier modes while 
$n(t)$ only has a finite number of modes.
Nevertheless, in numerical modeling, we can 
use $n(t)$ to approximate $\dot{{R}}(t)$.

Let us define
\begin{equation}
    \beta_t := \int_0^t n(t') dt' 
\end{equation}
and
\begin{equation}
    \Delta \beta_j = \beta_{t_j+\Delta t} - \beta_{t_j},
\end{equation}
and consider $ \Delta \beta_j$ as an approximation to $ \Delta {R}_j$ \textcolor{black}{(recall that $ \Delta {R}_j= R_{t_j+\Delta t}-R_{t_j}$ is the increment of the stochastic process $R$)}.
For the above-defined $n(t)$, we find
\begin{gather} \notag
\mathbb{E}\big[ \left(\Delta\beta_j\right)^2 \big]= \\
\mathbb{E}\left[\left\{
  \int_{t_j}^{t_j+\Delta t} \frac{1}{\sqrt{N_f\Delta t}}
  \left(
        C(\omega_0)\frac{b_0}{\sqrt{2}} 
    +  \sum_{m=1}^{N_f} C(\omega_m) 
       \big[ a_m\sin(\omega_m t) + 
       b_m\cos(\omega_m t)\big]\right)dt\right\}^2 \right] \notag
\end{gather}
For small $\Delta t,$ we can approximate $\mathbb{E}\big[ \left(\Delta\beta_j\right)^2 \big]$ as 
$$ 
\mathbb{E}\left[\left\{ \frac{1}{\sqrt{N_f\Delta t}}
  \left(
        C(\omega_0)\frac{b_0}{\sqrt{2}}
    +  \sum_{m=1}^{N_f} C(\omega_m) 
       \big[ a_m\sin(\omega_m t_j) + b_m\cos(\omega_m t_j)\big]
  \right)\Delta t\right\}^2
\right].
$$
Taking into account the independence among the coefficients $a_m$ and $b_m$, we have
\begin{eqnarray}
\mathbb{E}\big[ \left(\Delta\beta_j\right)^2 \big] 
&=&
\frac{1}{N_f\Delta t} 
\mathbb{E}\left[\,\left(
 C(\omega_0)\frac{b_0}{\sqrt{2}}\Delta t
 \right)^2  \right] \nonumber\\
&& 
  + \frac{1}{N_f\Delta t} \sum_{m=1}^{N_f} \mathbb{E}\left[ 
   C(\omega_m)^2 a_m^2 \sin^2(\omega_m t_j) (\Delta t)^2 \right] \nonumber \\
&&
+ \frac{1}{N_f\Delta t} \sum_{m=1}^{N_f} \mathbb{E}\left[
  C(\omega_m)^2 b_m^2 \cos^2(\omega_m t_j) (\Delta t)^2\right]
\end{eqnarray}
Also note that per construction, we have
\begin{equation}
    \mathbb{E}\left[ a_m^2 \right] \equiv 1 \,,\quad  \mathbb{E}\left[ b_m^2 \right] \equiv 1 \,,
\end{equation}
for any $m = 0, 1, \dots, N_f$.
Therefore 
\begin{gather}
\mathbb{E}\big[ \left(\Delta\beta_j\right)^2 \big] 
=
\frac{(\Delta t)^2}{N_f\Delta t} \left\{
 \frac{C(\omega_0)^2}{2} 
+ \sum_{m=1}^{N_f} C(\omega_m)^2\big[\sin^2(\omega_m t_j) + \cos^2(\omega_m t_j)\big]
\right\}
\label{eq:colored_noise_variance_2}\\
=
\Delta t \,\frac{1}{N_f} \left\{
 \frac{C(\omega_0)^2}{2}
+ \sum_{m=1}^{N_f} C(\omega_m)^2
\right\}.
\label{eq:colored_noise_variance_3}
\end{gather}
\textcolor{black}{The expression for $\mathbb{E}\big[ \left(\Delta\beta_j\right)^2 \big]$ can be used to obtain the {\it integral} form of the generalized It\^o correction, namely
\begin{gather}
		\lim_{\Delta t\rightarrow 0}\sum_j 
        \left(\frac 12 \frac{dg(u)}{du} g[u]\right) \!\!\bigg\vert_{t_j}\!\!
		\mathbb{E}\big[ \left(\Delta{R}_j\right)^2 \big]\\
		=
		\lim_{\Delta t\rightarrow 0}\sum_j 
        \left(\frac 12 \frac{dg(u)}{du} g[u]\right) \!\!\bigg\vert_{t_j}\!\! 
        \frac{1}{N_f} \left\{
 \frac{C(\omega_0)^2}{2}
+ \sum_{m=1}^{N_f} C(\omega_m)^2
\right\} \Delta t.
\end{gather}
Thus, the {\it differential} form of the generalized It\^o correction reads
\begin{eqnarray}
I = \frac 12 g(u)\frac{dg}{du}\,
\lim_{N_f \rightarrow \infty}\frac{1}{N_f}\! \left\{
 \frac{C(\omega_0)^2}{2} 
+ \sum_{m=1}^{N_f} C(\omega_m)^2
\right\}. \label{ito_correction_1}
\end{eqnarray}
The reason we keep in the final formula the limit $N_f \rightarrow \infty$ is that $N_f$ depends on the length of the subinterval $\Delta t.$ So, when $\Delta \rightarrow 0,$ the number of frequencies $N_f \rightarrow \infty.$ We have provided in the main text a short discussion on the behavior of this limit.} 

\begin{figure}[ht]
  \centering
\includegraphics[width = 0.9\textwidth]{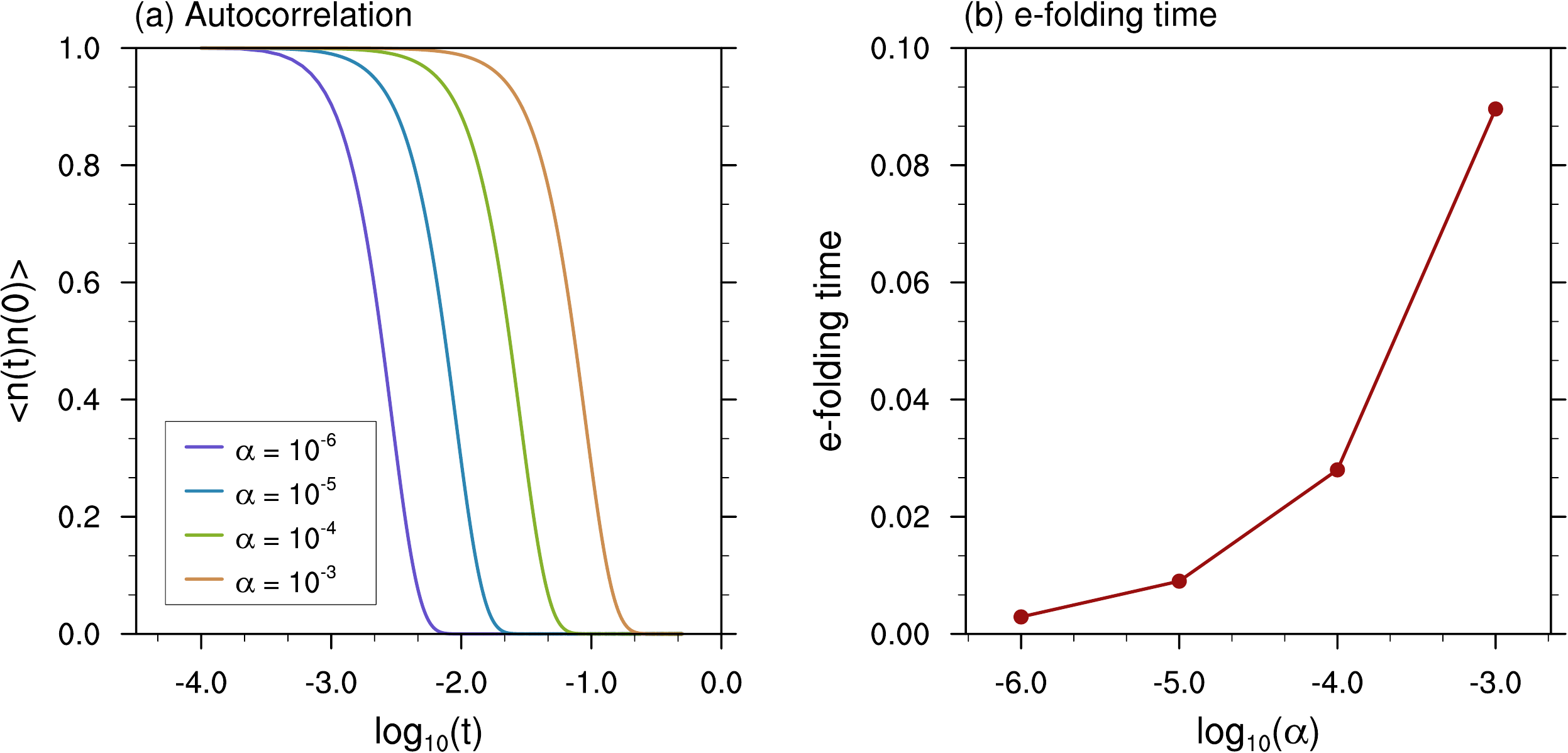}
  \caption{a) Autocorrelation of the noise process $n(t)$ for different values of $\alpha.$ b) The e-folding time of the noise process $n(t)$ for different values of $\alpha.$}
  \label{figure_colored_noise_statistics}
\end{figure}

\textcolor{black}{We have computed the autocorrelation and the e-folding time of the noise process $n(t)$ for different values of the parameter $\alpha$ (see Fig. \ref{figure_colored_noise_statistics}).}

\section*{Appendix B}
  \label{sec:spatial_analytical}


%

%
Let us define
\begin{eqnarray}
\mathbf{F} &=& (F_{-N_x}, F_{-N_x+1},\dots, F_{N_x-1}, F_{N_x})^T \\
\pmb{\mathsf{H}} &=& i\,\textrm{Diag}\{-N_x, -N_x+1,\cdots,N_x-1 , N_x\}
\end{eqnarray}
and
\begin{equation*}
\pmb{\mathsf{D}} = \left[ \begin{array}{cccc}
icN_x-\mu N_x^2 & \frac{i(N_x-1)\epsilon}{4} & 0 & 0\\
\frac{iN_x\epsilon}{4} & ic(N_x-1)-\mu (N_x-1)^2 & \frac{i(N_x-2)\epsilon}{4} & 0  \\
0 & \frac{i(N_x-1)\epsilon}{4} &  ic(N_x-2)-\mu (N_x-2)^2 & \frac{i(N_x-3)\epsilon}{4} \\
\cdots&\\
0 & \cdots  \\
0 & \cdots \\
0 & \cdots  
\end{array} \right.
\end{equation*}
\begin{equation*}
\left. \begin{array}{cccc} 0& 0  & \cdots & 0 \\
0 & 0  &\cdots &0\\
 \dots\\
 -\frac{i(N_x-3)\epsilon}{4} & -ic(N_x-2)-\mu(N_x-2)^2 & -\frac{i(N_x-1)\epsilon}{4} & 0 \\
  0 & -\frac{i(N_x-2)\epsilon}{4} &  -ic(N_x-1)-\mu(N_x-1)^2 & -\frac{iN_x\epsilon}{4} \\
 0 & 0 & -\frac{i(N_x-1)\epsilon}{4} &  -icN_x-\mu N_x^2
 \end{array}\right] 
\end{equation*}
This notation allows us to
write Eqs.~\eqref{eq:n}--\eqref{eq:-N_x} as Eq.~\eqref{eq:Fourier_ODE_nonconstant}
and the analytical solution as Eq.~\eqref{eq:analytic_solution_nonconstant}.

\bibliographystyle{plain}
\bibliography{reference}

\end{document}